\documentclass[12pt]{article}
\usepackage{latexsym}
\usepackage{graphicx}
\usepackage[dvips]{epsfig}
\usepackage{amsfonts}
\usepackage{amssymb}

\newtheorem{prob}{Problem}
\newcommand{\hole}{{\mkern2mu|\mkern2mu}}
\newcommand{\Ga}{{\mit \Gamma}}

\newcommand{\CP}{{\cal P}}
\newcommand{\CQ}{{\cal Q}}

\newcommand{\CG}{{\cal G}}
\newcommand{\CF}{{\cal F}}
\newcommand{\BE}{\mathbb{E}}  
  
\newcommand{\BS}{\mathbb{S}}

\newcommand{\BT}{\mathbb{T}} 
\newcommand{\BP}{\mathbb{P}} 
\newcommand{\BC}{\mathbb{C}} 
\newcommand{\eqref}[1]{(\ref{#1})}

\begin{document}

\title{Problems on Polytopes, Their Groups, and Realizations}

\author{Egon Schulte\thanks{Supported by NSA-grant H98230-05-1-0027} 
\hskip.05in and 
Asia Ivi\'{c} Weiss\thanks{Supported by NSERC of Canada Grant \#8857}}

\maketitle

\begin{abstract}
\noindent
The paper gives a collection of open problems on abstract polytopes that were either presented at 
the \emph{Polytopes Day in Calgary} or motivated by discussions at the preceding \emph{Workshop on Convex and Abstract Polytopes\/} at the Banff International Research Station in May~2006.
\end{abstract}

\section{Introduction}

The rapid development of polytope theory in the past thirty years has resulted in a rich theory featuring an attractive interplay of several mathematical disciplines. The breadth of the talks at the \emph{Workshop on Convex and Abstract Polytopes\/} and the subsequent \emph{Polytopes Day in Calgary} that we organized jointly with Ted Bisztriczky at the Banff International Research Station (BIRS) on May 19-21, 2005 and the University of Calgary on May 22, 2005, respectively, gave evidence that polytope theory is very much alive and is the unifying theme of a lot of research activity. The Workshop provided a much desired opportunity to share recent developments and emerging directions on geometric, combinatorial, and abstract aspects of polytope theory.  It is noteworthy that the last major meeting on convex and abstract polytopes was the NATO Advanced Study Institute on \emph{Polytopes -- Abstract, Convex and Computational\/} in 1993 in Scarborough, Ontario (see \cite{bis}). 

For abstract polytopes, the invited lectures and talks focused on polytopes with various degrees of combinatorial or geometric symmetry (regular, chiral, or equivelar polytopes, and their geometric realization theory), as well as the structure of their symmetry groups or automorphism groups (reflection groups, Coxeter groups, and C-groups, and their representation theory).  The present paper surveys open research problems on abstract polytopes that were presented at the Workshop (primarily at the problem session). A significant number of these problems have been addressed in detail elsewhere in the literature, notably 
in~\cite{arp}. For the remaining problems we provide some background information when available, but due to space limitations we cannot give a comprehensive account in all cases.

In Section~\ref{bano} we review basic notions and concepts, and then in the subsequent sections explore some of the most important problems in the field. It is natural to group these problems under the headings and subheadings provided in Sections~\ref{polmap} to \ref{realiz}. There are, however, a number of interesting research problems that do not fall into these sections; they are presented in Section~\ref{othpro}.  

\section{Basic Notions}
\label{bano}

For the general background on abstract polytopes we refer the reader to the monograph~\cite[Chapters~2, 3]{arp}. Here we just review some basic terminology.

An (\emph{abstract\/}) \emph{polytope of rank\/} $n$, or an \emph{$n$-polytope\/}, is a partially ordered set $\CP$ with a strictly monotone rank function with range $\{-1,0, \ldots, n\}$. The elements of $\CP$ with rank $j$ are the \emph{$j$-faces\/} of $\cal P$, and the faces of ranks $0$, $1$ or $n-1$ are also called \emph{vertices\/}, \emph{edges\/} or \emph{facets\/}, respectively. The maximal chains are the \emph{flags} of $\CP$ and contain exactly $n + 2$ faces, including a unique minimal face $F_{-1}$ (of rank $-1$) and a unique maximal face $F_n$ (of rank $n$). Two flags are called \emph{adjacent} ($i$-\emph{adjacent}) if they differ in just one face (just their $i$-face, respectively); then $\CP$ is \emph{strongly flag-connected}, meaning that, if $\Phi$ and $\Psi$ are two flags, then they can be joined by a sequence of successively adjacent flags $\Phi = \Phi_0,\Phi_1,\ldots,\Phi_k = \Psi$, each of which contains $\Phi \cap \Psi$. Furthermore, $\CP$ has the following homogeneity property:\  whenever $F \leq G$, with $F$ a $(j-1)$-face and $G$ a $(j+1)$-face for some $j$, then there are exactly two $j$-faces H with $F \leq H \leq G$. These conditions essentially say that $\CP$ shares many combinatorial properties with classical polytopes.

For any two faces $F$ and $G$ with $F \leq G$ we call $G/F :=\{ H \mid F \leq H \leq G \}$ a \emph{section\/} of $\CP$; this is a polytope (of the appropriate rank) in its own right. In general, there is little possibility of confusion if we identify a face $F$ and the section $F/F_{-1}$. If $F$ is a face, then $F_{n}/F$ is said to be the \emph{co-face at\/} $F$, or the \emph{vertex-figure at\/} at $F$ if $F$ is a vertex.

A polytope $\CP$ is \emph{regular\/} if its combinatorial automorphism group $\Ga(\CP)$ is transitive on the flags of $\CP$. The group of a regular $n$-polytope $\CP$ is generated by involutions $\rho_0,\ldots,\rho_{n-1}$, where $\rho_i$ maps a fixed, or \emph{base\/}, flag $\Phi$ to the flag $\Phi^i$, $i$-adjacent to $\Phi$. These \emph{distinguished generators\/} satisfy (at least) the standard Coxeter-type relations 
\begin{equation}
\label{standrel}
(\rho_i \rho_j)^{p_{ij}} = \epsilon \;\; \textrm{ for } i,j=0, \ldots,n-1, 
\end{equation}
where $p_{ii}=1$, $p_{ji} = p_{ij} =: p_{i+1}$ if $j=i+1$, and $p_{ij}=2$ otherwise; thus the underlying Coxeter diagram is a string diagram. The numbers $p_{j} := p_{j-1,j}$ ($j = 1,\ldots,n-1$) determine the (\emph{Schl\"afli}) \emph{type} $\{p_{1},\ldots,p_{n-1}\}$ of $\CP$. The distinguished generators also satisfy the {\em intersection condition\/},
\begin{equation}
\label{intprop}
\langle \rho_i \mid i \in I \rangle \cap \langle \rho_i \mid i \in J \rangle 
= \langle \rho_i \mid i \in {I \cap J} \rangle 
\qquad (I,J \subseteq \{0,1,\ldots,n-1\}). 
\end{equation}

Groups $\Ga = \langle \rho_0, \ldots, \rho_{n-1} \rangle$ whose generators satisfy \eqref{standrel} and \eqref{intprop}, are called \emph{string C-groups\/}; here, the ``C" stands for ``Coxeter", though not every C-group is a Coxeter group. These string C-groups are precisely the automorphism groups of regular polytopes. In this context it is often convenient to adopt the viewpoint that a regular polytope $\CP$ is to be identified with its string C-group $\Ga := \Ga(\CP)$ (see \cite[Theorem~2E11]{arp}). The \emph{$j$-faces} of $\CP$ then are the right cosets $\Ga_j \varphi$ of the \emph{distinguished subgroup} $\Ga_j := \langle \rho_i \mid i \neq j \rangle$ for each $j=0,\ldots,n-1$, and two faces are incident just when they intersect (as cosets), that is, $\Ga_j \varphi \leq \Ga_k\psi$ if and only if $\Ga_j\varphi \cap \Ga_k\psi \neq \emptyset$ and $j \leq k$. Formally, we must also adjoin two copies of $\Ga$ itself, as the (unique) $(-1)$-face and $n$-face of $\CP$. If the regular polytope is determined just by the $p_j$, then we have the \emph{universal} regular polytope (with that Schl\"afli type), for which we use the same symbol $\{p_1,\ldots,p_{n-1}\}$. We write $[p_1,\ldots,p_{n-1}]$ for the corresponding \emph{Coxeter} group. Generally, however, the group $\Ga$ will satisfy additional relations as well.

\begin{figure}[hbt]
\begin{center}
\begin{picture}(150,200)
\put(-95,0)
{\includegraphics[width=5in]{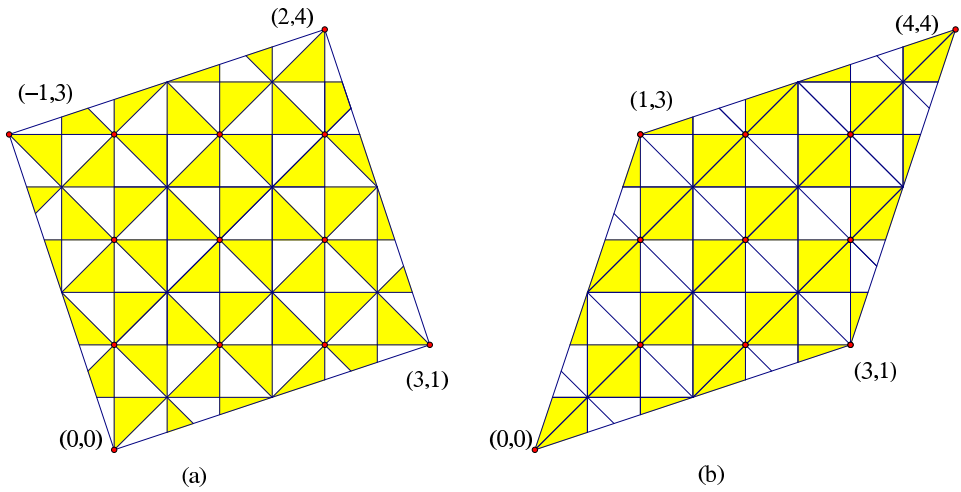}}
\end{picture}
\caption{\em Two torus maps of type $\{4,4\}$}
\label{figone}
\end{center}
\end{figure}

A polytope $\CP$ (of rank $n \geq 3$) is {\em chiral\/} if its group $\Ga(\CP)$ has two orbits on the flags, such that adjacent flags belong to distinct orbits (see \cite {arp,swp}). Figure~\ref{figone}(a) shows the chiral torus map $\{4,4\}_{(3,1)}$ obtained from the square tessellation in the plane by factoring out the lattice spanned by the vectors $(3,1)$ and $(-1,3)$ (see \cite{cm}). Chiral polytopes form an important class of nearly regular polytopes. Intuitively, they have complete (combinatorial) rotational symmetry, but not the full symmetry by (combinatorial) reflection, as highlighted by the following local characterization of chirality.

A polytope $\CP$ is chiral if and only if $\CP$ is not regular, but for some base flag $\Phi = \{F_{-1}, F_0, \ldots, F_n \}$ of $\CP$ there exist automorphisms $\sigma_1,
\ldots, \sigma_{n-1}$ of $\CP$ such that $\sigma_i$ fixes all faces in $\Phi \setminus \{
F_{i-1}, F_{i} \}$ and cyclically permutes consecutive $i$-faces of $\CP$ in the section
$F_{i+1} / F_{i-2}$ (of rank $2$). For a chiral polytope $\CP$, we can choose the orientation of the $\sigma_{i}$ for $i = 1,\ldots,n-1$ in such a way that, if $F_{i}'$ denotes the $i$-face of $\CP$ with $F_{i-1} < F_{i}' < F_{i+1}$ and $F_{i}' \neq F_{i}$, then $F_{i}\sigma_{i} = F_{i}'$, and thus $F_{i-1}'\sigma_{i} = F_{i-1}$.  Then $\sigma_{1},\ldots,\sigma_{n-1}$ generate the group $\Ga(\CP)$ and satisfy (at least) the relations
\begin{equation}
\label{relchir}
\begin{array}{rcl}
\sigma_i^{p_i}  & = & \epsilon\;\, \textrm{ for } 1\leq i \leq n-1, \\[.03in]
(\sigma_i \sigma_{i+1}\dots \sigma_j)^2 & = & \epsilon\;\, \textrm{ for } 1\leq i < j \leq n-1,
\end{array} 
\end{equation}
where again $\{p_1,\ldots,p_{n-1}\}$ is the {\em type} of $\CP$.  The \emph{distinguished generators\/} $\sigma_1,\ldots, \sigma_{n-1}$ of $\Ga(\CP)$ satisfy an intersection condition which (especially for $n>4$) is more complicated than for regular polytopes (see \cite{swp}). When $n=4$ this condition takes the form
\begin{equation}
\label{intchir4}
\begin{array}{c}
\langle\sigma_1\rangle \cap \langle\sigma_2\rangle = \{\epsilon\} = 
\langle \sigma_2 \rangle \cap \langle \sigma_3 \rangle,\quad
\langle \sigma_1, \sigma_2 \rangle \cap \langle  \sigma_2, \sigma_3 \rangle 
= \langle \sigma_2 \rangle
\end{array}
\end{equation}
while for $n=3$ we only have 
\begin{equation}
\label{intchir3}
\langle \sigma_1  \rangle \cap \langle \sigma_2 \rangle = \{\epsilon\}.
\end{equation}

For a regular polytope $\CP$ with group $\Ga(\CP) = \langle \rho_{0},\rho_{1},\ldots,\rho_{n-1} \rangle$, define $\sigma_{i} := \rho_{i-1}\rho_{i}$ for $i=1,\ldots,n-1$. Then $\sigma_{1},\ldots,\sigma_{n-1}$ generate the (\emph{combinatorial\/}) \emph{rotation subgroup\/} $\Ga^{+}(\CP)$ of $\Ga(\CP)$ (of index at most $2$) and have properties similar to those of the distinguished generators for chiral polytopes. We call a regular polytope $\CP$ \emph{directly regular\/} if $\Ga^{+}(\CP)$ has index $2$ in $\Ga(\CP)$ 

If $\Ga = \langle \sigma_1, \ldots, \sigma_{n-1} \rangle$ is a group whose generators satisfy  relations as in \eqref{relchir} as well as the respective intersection condition (that is, \eqref{intchir4} if $n=4$, or \eqref{intchir3} if $n=3$), then there exists a chiral polytope $\CP$ with $\Ga(\CP) \cong \Ga$ or a directly regular polytope $\CP$ with $\Ga^{+}({\CP}) \cong \Ga$. Moreover, $\CP$ is directly regular if and only if there is an involutory group automorphism $\alpha: \Ga \rightarrow\Ga$ such that $(\sigma_1)\alpha = \sigma_1^{-1}$, $(\sigma_2)\alpha =\sigma_1^{2} \sigma_2$
and $(\sigma_i)\alpha = \sigma_i$ for $i=3,\ldots,n-1$.

Each isomorphism class of a chiral polytope occurs in two \emph{enantiomorphic forms}, each determined by a choice (of orbit) of base flag (under the group). Thus the two enantiomorphic forms of a chiral polytope $\CP$ are represented by a pair of (non-equivelant) adjacent base flags for $\CP$, or, alternatively, two (non-equivalent) systems of generators for the group $\Ga(\CP)$. If the generators of $\Ga(\CP)$ associated with the base flag $\Phi$ are given by $\sigma_1,\ldots,\sigma_{n-1}$, then those associated with the $0$-adjacent flag $\Phi^{0}$ are $\sigma_1^{-1}, \sigma_1^{2}\sigma_{2}, \sigma_3,\ldots,\sigma_{n-1}$. See \cite {swp} for more details. (For directly regular polytopes, the respective enantiomorphic forms still can be defined but can naturally be identified.)

\section{Polyhedra and maps}
\label{polmap}

Polytopes of rank $3$ are also called \emph{polyhedra\/} (see \cite[Section 7B]{arp}). Each polyhedron $\CP$ yields a map (tessellation) on a closed surface; and vice versa, most maps on surfaces arise in this way, in which case they are called \emph{polytopal}. 

We begin here with rank $3$. Even though most of the problems can be formulated for polytopes of any rank, we deliberately concentrate on polyhedra, as, although difficult, the problems so posed have an intuitive appeal. A considerable number of papers have been published on regular and chiral maps; for classification results by genus see, for example, \cite{confmaps,cm,ga,sh} and more recently \cite{conder}. Here we touch on three specific topics related to maps, which can serve to illustrate how some easily stated problems in geometry can present major challenges.

\bigskip
\noindent
{\bf Self-duality\/}
\vskip.02in

Recall that a polyhedron $\CP$ is called \emph{self-dual} if and only if it is isomorphic to its dual; in this case $\CP$ possesses a \emph{duality}, that is, an incidence-reversing bijection. A duality of period $2$ is called a \emph{polarity}.

\begin{prob}
\label{selftor} 
Classify the self-dual toroidal polyhedra.
\end{prob}

Similar such problems can of course be formulated for polyhedra of any genus. The spherical self-dual polyhedra have been classified independently in \cite{arri} and \cite{serv}. 

Call a polyhedron $\CP$ \emph{equivelar} of type $\{p,q\}$ if all its $2$-faces are $p$-gons and all its vertex-figures are $q$-gons.

\begin{prob}
\label{equipol} 
Do equivelar self-dual polyhedra possess a polarity?
\end{prob}

The answer is affirmative for regular and chiral polyhedra (see \cite{self-duality}) but is not known even for the fully transitive polyhedra discussed below. In general, self-dual polyhedra need not admit a polarity, as was shown in, for example, \cite{jen}, answering a question posed in \cite{grshi}.

\bigskip
\noindent
{\bf Transitivity\/}
\vskip.02in

While in the classical theory of convex polyhedra equivelarity is sufficient to characterize the regular polyhedra, this is certainly not true for abstract polyhedra, as already the existence of chiral polyhedra shows. The most symmetric equivelar polyhedra have automorphism groups that act transitively on the faces of each rank. We shall call such polyhedra \emph{fully transitive\/}. Every regular or chiral polyhedron is fully transitive.

\begin{prob}
\label{ftclas} 
Classify fully transitive polyhedra on orientable or non-orientable surfaces of small genus.
\end{prob}

This has been done for the torus (see \cite{hubtor}). In addition to the regular and chiral torus maps, there is just one infinite family of fully transitive toroidal maps, of type $\{4,4\}$, one member of which is depicted in Figure~\ref{figone}(b) and is obtained from the square tessellation by factoring out the lattice spanned by the vectors $(3,1)$ and $(1,3)$; as for chiral maps (such as in Figure~\ref{figone}(a)), there are just two flag orbits, but now every two $1$-adjacent flags are in the same orbit. It is well-known that there are no chiral maps on orientable surfaces of genus $g=2$, $3$, $4$, $5$ or $6$ (and none on non-orientable surfaces), so it is natural to pose the following

\begin{prob}
\label{smallgen} 
What is the smallest integer $g$, $g\geq 2$, for which there exists a non-regular fully transitive polyhedron on an orientable surface of genus $g$?
\end{prob}

The groups of regular polyhedra are characterized by equations \eqref{standrel} and \eqref{intprop}, and those of chiral polyhedra by equations \eqref{relchir} and \eqref{intchir3}.

\begin{prob}
\label{chargr} 
Can one similarly characterize the automorphism groups for fully-transitive polyhedra (which are neither regular nor chiral)?
\end{prob}

If Problem~\ref{chargr} seems too difficult, one may attempt to find a characterization for self-dual polyhedra only.

The analogous problems for edge-transitive polyhedra are obviously even more challenging. Partial results in this direction have recently been obtained in \cite{grwat,orb,situ}. 

\begin{prob}
\label{clasedge} 
Classify edge-transitive polyhedra on orientable or non-orientable surfaces of small genus.
\end{prob}

\medskip
\noindent
{\bf Zigzags and holes\/}
\vskip.02in

A \emph{$j$-zigzag\/} of a polyhedron $\CP$ is an edge-path on the underlying surface which leaves a vertex by the $j$-th edge from which it entered, alternating the sense (in some local orientation). Similarly, a \emph{$j$-hole\/} of $\CP$ is an edge-path on the surface which leaves a vertex by the $j$-th edge from which it entered, always in the same sense, keeping to the left (say); we simply refer to it as a \emph{hole\/} if $j=2$. For example, the $1$-zigzags are precisely the \emph{Petrie polygons\/} of $\CP$, and the $1$-holes are simply the $2$-faces of $\cal P$ (see \cite{cm}). If $\CP$ is regular, the lengths of a $j$-zigzag and $j$-hole are given by the periods of $\rho_0 (\rho_1 \rho_2)^j$ or $\rho_0\rho_1(\rho_2\rho_1)^{j-1}$ in $\Ga(\CP)$, respectively. Similarly, if $\CP$ is chiral, the length of a $j$-hole is the period of $\sigma_1 \sigma_2^{1-j}$ in $\Gamma({\cal P})$. On the other hand, a chiral polyhedron $\CP$ is orientable, so a $j$-zigzag must necessarily have even length (as orientation is reversed at each step along it); this length is twice the period of $\sigma_{1}\sigma_{2}^{-j}\sigma_{1}^{-1}\sigma_{2}^{j}$ or, alternatively, its conjugate $\sigma_{2}^{j}\sigma_{1}\sigma_{2}^{-j}\sigma_{1}^{-1}$ in $\Ga(\CP)$, so that ``right" and ``left" $j$-zigzags of a chiral polyhedron always have the same lengths.

If a regular polyhedron of type $\{p, q\}$ is completely determined by the lengths $h_j$ of its $j$-holes, for $2 \leq j \leq  k:= \lfloor q/2 \rfloor$, and lengths $r_j$ of its $j$-zigzags, for $1 \leq j \leq  k:= \lfloor q/2 \rfloor$, then we denote it by
\begin{equation}
\label{hozig}
{\cal P} = \{p,q \hole h_2, \ldots, h_k \}_{r_1, \ldots,r_k}, 
\end{equation}
with the convention that any unnecessary $h_j$ or $k_j$ (that is, one that need not be 
specified) is replaced by a $\cdot$, with those at the end of the sequence omitted (see \cite[Section 7B]{arp}). More generally, we say that a regular or chiral polyhedron is of 
{\em type\/} $\{p,q \hole h_2, \ldots, h_k \}_{r_1,\ldots,r_k}$ if its $j$-holes and $j$-zigzags have lengths $h_j$ or $r_j$, respectively, with $j$ as above. 

\begin{prob}
\label{zighole} 
To what extent is a regular or chiral polyhedron of type $\{p,q\}$ determined by the lengths $h_j$ of its $j$-holes and the lengths $r_j$ of its $j$-zigzags, for $1 \leq j \leq  k:= \lfloor q/2 \rfloor$? In other words, what can be said about the quotients of the regular polyhedron 
$\{p,q \hole h_2, \ldots, h_k \}_{r_1, \ldots,r_k}$ defined in \eqref{hozig}? In particular, when is the latter finite? 
\end{prob}

\section{The Amalgamation Problem}
\label{amalpro}

A main thrust in the theory of abstract polytopes, as well as it was in the classical theory, is that of the amalgamation of polytopes of lower rank (see \cite[Chapter 4]{arp}). In a gist, the amalgamation problem can be described as follows. We concentrate on polytopes with high symmetry, although the problem makes sense in a more general context.

\begin{prob}
\label{am} 
Given regular (or chiral) polytopes $\CP_1$ and $\CP_2$ of rank $n$, does there exist a regular (or chiral) polytope $\CP$ of rank $n+1$ with facets isomorphic to $\CP_1$ and vertex-figures isomorphic to $\CP_2$?
\end{prob}

Here the vertex-figures of $\CP_1$ must be isomorphic to the facets of $\CP_2$. Even then, the general answer is negative; for example, there is no regular polytope with toroidal facets $\{4,4\}_{(3,0)}$ and hemi-cubes $\{4,3\}_3$ as vertex-figures (see \cite[p. 367]{arp}). Thus we arrive at the following refined version of the amalgamation problem. 

\begin{prob} 
\label{amre}
Describe or characterize classes of regular (or chiral) $n$-polytopes $\CP_1$ and $\CP_2$ with an affirmative answer to Problem~\ref{am}.
\end{prob}

If any regular or chiral $(n+1)$-polytope $\CP$ with regular facets $\CP_1$ and regular vertex-figures $\CP_2$ exists, then there is a universal such polytope, denoted $\{\CP_1,\CP_2\}$, which is regular. 

\begin{prob}
\label{reguni}
For regular $n$-polytopes $\CP_1$ and $\CP_2$ with an affirmative answer to Problem~\ref{am}, determine the group of the universal $(n+1)$-polytope $\{\CP_1,\CP_2\}$ and give conditions under which it is finite. 
\end{prob}

There has been much interest in these problems in the context of the topological classification of regular polytopes (see \cite[Chs. 4, 10--12]{arp}). However, in their full generality the questions are still open. 

If the facets or vertex-figures are chiral, then the situation is more subtle (see \cite{swp}). Here the outcome depends on the choice of enantiomorphic form for $\CP_1$ and $\CP_2$. This distinction is largely irrelevant if at least one of $\CP_1$ and $\CP_2$ is directly regular, but it is essential if both $\CP_1$ and $\CP_2$ are chiral. It still is true that, if there exists any chiral $(n+1)$-polytope $\CP$ with preassigned enantiomorphic form of chiral or (necessarily, directly) regular facets $\CP_1$ and vertex-figures $\CP_2$ (where not both are directly regular), then there is also a universal such polytope, denoted by $\{\CP_1,\CP_2\}^{ch}$; here, by slight abuse of notation, we have let $\CP_1$ and $\CP_2$ denote the enantiomorphic forms rather than the polytopes themselves. Thus, given isomorphism classes of chiral or regular facets and vertex-figures, there are generally four amalgamation problems, but modulo enantiomorphism they reduce to only two problems (that is, the solutions for one are enantiomorphic to solutions for another), or to only one problem if $\CP_1$ or $\CP_2$ is regular. For example, the chiral universal $4$-polytopes $\{\{4,4\}_{(1,3)},\{4,4\}_{(1,3)}\}$ and $\{\{4,4\}_{(3,1)},\{4,4\}_{(1,3)}\}$ are not isomorphic; their groups are of orders $2000$ and $960$, respectively (see \cite[p. 132]{self-duality}).

\begin{prob}
\label{reguni}
As before, let $\CP_{1},\CP_{2}$ denote enantiomorphic forms of two chiral or regular $n$-polytopes which are not both regular.  For any such pair $\CP_{1},\CP_{2}$ with an affirmative answer to Problem~\ref{am}, determine the group of the universal $(n+1)$-polytope $\{\CP_1,\CP_2\}^{ch}$ and give conditions under which it is finite. 
\end{prob}

The next two problems deal with a particularly interesting case.

\begin{prob}
\label{resphvert}
Determine existence and finiteness for the universal regular $(n+1)$-polytopes $\{\CP_1,\CP_2\}$, where $\CP_1$ is a finite (abstract) regular $n$-polytope and $\CP_2$ is a convex regular $n$-polytope (regular spherical tessellation on $\BS^{n-1}$).
\end{prob}

Problem~\ref{resphvert} has already been solved when the vertex-figure $\CP_2$ is an $n$-cross-polytope (see \cite[Thm. 8E10]{arp}); in this case the universal $(n+1)$-polytope is finite if and only if $\CP_1$ is \emph{neighborly}, meaning here that any two vertices of $\CP_1$ are joined by an edge.

\begin{prob}
\label{chsphvert}
Determine existence and finiteness for the universal chiral $(n+1)$-polytopes  $\{\CP_1,\CP_2\}^{ch}$, where $\CP_1$ is (an enantiomorphic form of) a finite (abstract) chiral $n$-polytope and $\CP_2$ is a convex regular $n$-polytope (spherical regular tessellation on $\BS^{n-1}$).
\end{prob}

Both Problems~\ref{resphvert} and \ref{chsphvert} are open even in the special case when $\CP_2$ is a $d$-simplex. Even though this seems difficult enough, a much more challenging problem might be the case when the vertex-figure $\CP_2$ is allowed to be a regular $n$-star-polytope (see \cite[Sect. 7D]{arp}).

The amalgamation problems discussed here ask for the existence and properties of regular or chiral polytopes with preassigned regular or chiral facets and vertex-figures. If only one kind, the facets (say), are prescribed, we arrive at \emph{extension problems\/} for regular or chiral polytopes. Here the basic problem is to find and describe polytopes with a given isomorphism type of facet. Several general extension results have been established in the literature (see, for example, \cite[Sections 4D, 8C,D]{arp}), but most concentrate on the regular case. Recent progress on extensions of regular polytopes with prescribed last entry of Schl\"afli symbol has been obtained in~\cite{peco}. Interestingly, symmetric or alternating groups, or products of such groups, often occur as automorphism groups of the extending polytopes (for example, 
in~\cite{pel}, polyhedra with alternating groups are constructed).

Many of the above problems are contingent upon the following open problem, for which we conjecture an affirmative answer.

\begin{prob}
\label{chirank}
Do finite chiral polytopes of rank $n\geq 6$ exist?
\end{prob}

In rank $5$ two such polytopes have recently been constructed in \cite{cohupi}. Known examples in rank $4$ are numerous.

The last problem in this section is a long-standing open problem about convex $4$-polytopes posed in \cite{pesh} (see also \cite{kal,esnon}). 

\begin{prob}
\label{icos}
Does there exist a convex $4$-polytope all of whose facets are combinatorially isomorphic to the icosahedron? 
\end{prob}

A convex $d$-polytope $P$ is called a \emph{facet\/} if there exists a convex $(d+1)$-polytope all of whose facets are combinatorially isomorphic to $P$; otherwise, $P$ is a \emph{nonfacet}.  Problem~\ref{icos} asks whether or not the icosahedron is a facet or nonfacet.  For a detailed study of facets and nonfacets of convex polytopes (and tiles and nontiles of Euclidean spaces) see \cite{pesh,esnon}.

\section{Topological classification}
\label{topclass}

In contrast to the traditional theory where a convex polytope is locally and globally spherical, it is a very subtle problem to define the topological type of an abstract polytope (see \cite[Chapter 6]{arp}). In fact, this cannot be done unambiguously, except in certain cases. 

Call an abstract $n$-polytope $\CP$ (\emph{globally}) \emph{spherical}, \emph{projective}, or \emph{toroidal\/} if it is isomorphic to the face-set of a locally finite face-to-face tessellation on the $(n-1)$-sphere $\BS^{n-1}$, projective $(n-1)$-space $\BP^{n-1}$, or the $(n-1)$-torus $\BT^{n-1}$, respectively. The spherical regular $n$-polytopes are the convex regular $n$-polytopes (realized as tessellations on $\BS^{n-1}$), and the projective regular $n$-polytopes are the quotients of centrally symmetric convex regular $n$-polytopes by their central symmetry (see \cite{crp} and \cite[Section 6C]{arp}). The toroidal regular polytopes, or \emph{regular toroids}, of rank $n$ have been described in \cite{confmaps,cm} for $n=2$, and in \cite[Sections 6E,F]{arp} for $n\geq 3$. They are quotients of regular tessellations in Euclidean $(n-1)$-space $\BE^{n-1}$ by certain lattices, namely scaled copies of the hexagonal lattice in $\BE^2$, or the cubic lattice, the face-centered cubic lattice (root lattice $D_{n-1}$), or the body-centered cubic lattice (dual to $D_{n-1}$) in $\BE^{n-1}$ (the lattice $D_4$ in $\BE^4$ yields three families of regular toroids, up to duality). From the plane tessellations $\{3,6\}$, $\{6,3\}$ and $\{4,4\}$ we obtain the regular polyhedra $\{3,6\}_{(s,t)}$, $\{6,3\}_{(s,t)}$ or $\{4,4\}_{(s,t)}$ (with $t=0$ or $s=t$) on the $2$-torus, respectively. The cubical tessellation $\{4,3^{n-3},4\}$ in $\BE^{n-1}$, $n\geq 3$, yields the \emph{cubic} toroids $\{4,3^{n-3},4\}_{\bf s}$, where ${\bf s}:=(s^k,0^{n-1-k})$ with $s\geq 2$ and $k=1$, $2$ or $n-1$; when $n=3$ we obtain $\{4,4\}_{(s,0)}$ or $\{4,4\}_{(s,s)}$.  Finally, there are the regular toroids $\{3,3,4,3\}_{\bf s}$ and $\{3,4,3,3\}_{\bf s}$ of rank $5$, where ${\bf s}:=(s^k,0^{4-k})$ with $s\geq 2$ and $k=1$ or $2$, derived from the tessellations $\{3,3,4,3\}$ and $\{3,4,3,3\}$ in $\BE^4$. For all regular toroids elegant presentations for their automorphism groups are known.  By contrast, chiral polytopes cannot be spherical or projective for any rank $n$, and can also not be toroidal if $n>3$.  Moreover, the torus is the only compact Euclidean space-form which can admit a regular or chiral tessellation, and chirality can only occur on the $2$-torus. Except when $n=3$, little is known about regular or chiral $n$-polytopes of other topological types.

Given a topological type $X$, we say that an abstract $(n+1)$-polytope $\CP$ is \emph{locally of topological type $X$} if its facets and vertex-figures which are not spherical, are of topological type $X$. (There are variants of this definition, not adopted here, requiring the minimal sections of $\CP$ which are not spherical, to be of the specified topological type.)  In our applications, $X$ is $\BS^{n-1}$, $\BP^{n-1}$ or $\BT^{n-1}$ with $n\geq 3$, or $X$ is a $2$-dimensional surface and $n=3$. More generally, given two topological types $X_1$ and $X_2$ (of the kind just mentioned), an abstract $(n+1)$-polytope $\CP$ is said to be of \emph{local topological type $(X_{1},X_{2})$} if its facets are of topological type $X_1$ and its vertex-figures of topological type $X_2$.

\medskip
\noindent
{\bf Locally spherical or projective polytopes}
\smallskip

Every locally spherical regular or chiral $(n+1)$-polytope of type $\{p_{1},\ldots,p_{n}\}$ is a quotient of a (universal) regular tessellation $\{p_{1},\ldots,p_{n}\}$ in spherical, Euclidean or hyperbolic $n$-space (\cite[Section 6B]{arp}). In other words, locally spherical regular or chiral polytopes are well-understood up to taking quotients. Moreover, 
\[ \{p_{1},\ldots,p_{n}\} = \{\{p_{1},\ldots,p_{n-1}\}, \{p_{2},\ldots,p_{n}\}\}, \]
that is, the tessellation is also universal among all regular or chiral $(n+1)$-polytopes with spherical facets $\{p_{1},\ldots,p_{n-1}\}$ and vertex-figures $\{p_{2},\ldots,p_{n}\}$. In particular, finiteness occurs only when the tessellation itself is spherical.

The classification of the locally projective regular polytopes was recently completed in 
\cite{halp} (see also \cite[Section 14A]{arp}).  In rank $4$ (but not in rank $5$), all locally projective regular polytopes are finite. There are seventeen universal locally projective regular $4$-polytopes (including eight examples with a $2$ in their Schl\"afli symbol); amongst their $441$ quotients are a further four (non-universal) regular polytopes. A particularly interesting example is the universal $4$-polytope $\{\{5,3\},\{3,5\}_{5}\}$ (with dodecahedral facets and hemi-icosahedral vertex-figures), whose group is the direct product of two simple groups, the Janko group $J_1$ and $\rm{PSL}_{2}(19)$ (see \cite{halem}).

\medskip
\noindent
{\bf Locally toroidal regular polytopes}
\smallskip

Nearly thirty years ago, Gr\"unbaum~\cite{grgcd} posed the challenging problem (for $n=3$), as yet unsolved in full generality, of completely classifying the regular $(n+1)$-polytopes which are locally toroidal. Such polytopes can only exist for small ranks, namely $4$, $5$ or $6$.  Considerable progress has been made towards a complete classification (see \cite[Chapters 10--12]{arp}).  In this context, classification means the complete enumeration of all the locally toroidal, {\em finite universal\/} regular polytopes with prescribed types of facets and vertex-figures. In algebraic terms, this classification translates into the enumeration of certain quotients of hyperbolic string Coxeter groups defined in terms of generators and relations. At present the enumeration is complete in rank $5$, and nearly complete in rank $4$, while in rank $6$ there exist lists of finite universal regular polytopes strongly conjectured to be complete.  An up-to-date account on the state of the classification can be found in \cite{arp}. Tori and decompositions of spaces by tori are of such fundamental importance in topology that there is a good chance that progress on these problems would have impact on other areas in mathematics. 

In rank $4$, locally toroidal regular polytopes can have one of seven types (up to duality), namely 
\begin{equation}
\label{schla}
\{4,4,3\},\, \{4,4,4\},\, \{6,3,3\},\, \{6,3,4\},\, \{6,3,5\},\, \{6,3,6\} \mbox{ or } \{3,6,3\} . 
\end{equation} 
The enumeration is complete except in the cases $\{4,4,4\}$ and $\{3,6,3\}$. For example, for $\{4,4,4\}$ the classification amounts to the analysis of the universal regular polytopes $\{\{4,4\}_{(s,t)},\{4,4\}_{(u,v)}\}$; this involves answering the following questions: for which set of parameters $s,t,u,v$ do these polytopes exist, what are their groups, and when are they finite? For the type $\{4,4,4\}$ only the following choice of parameters has not been settled.

\begin{prob}
\label{fff}
Classify the universal regular $4$-polytopes $\{\{4,4\}_{(s,0)},{\{4,4\}}_{(u,0)}\}$ with $s,u\geq 3$, odd and distinct. It is conjectured (see \cite[p.376]{arp}) that these polytopes exist for all such $s$ and $u$, but are finite only if $(s,u) = (3,5)$ or $(5,3)$.
\end{prob}

For $\{3,6,3\}$, only partial results are known and involve sparse sequences of parameters (see \cite[Sections 11E,H]{arp}).

\begin{prob}
\label{tst}
Classify the universal regular $4$-polytopes $\{\{3,6\}_{(s,t)},{\{6,3\}}_{(u,v)}\}$. Here, 
$s\geq 2$, $t=0$ or $s=t\geq 1$, and $u\geq 2$, $v=0$ or $u=v\geq 1$. 
\end{prob}

In rank $5$, the enumeration is complete and only involves polytopes of type $\{3,4,3,4\}$ or $\{4,3,4,3\}$. 

However, in rank $6$, the problem is wide open and only partial results are known (see \cite[Sections 12C,D,E]{arp}). Up to duality, the respective universal regular $6$-polytopes are
\begin{equation} 
\label{sixpol}
\{\{3,3,3,4\},\{3,3,4,3\}_{\bf s}\},\;\,
\{\{3,3,4,3\}_{\bf s},\{3,4,3,3\}_{\bf t}\}\, \mbox{ or }\,
\{\{3,4,3,3\}_{\bf t},\{4,3,3,4\}_{\bf u}\},
\end{equation}
where ${\bf s}:=(s^k,0^{4-k})$, ${\bf t}:=(t^k,0^{4-k})$, ${\bf u}:=(u^l,0^{4-l})$ with $s,t,u\geq 2$ and $k=1,2$ or $l=1,2,4$. Now the facets and vertex-figures are regular toroids of rank $5$ (that is, tessellations on $4$-dimensional tori), except that the facets of the first kind are $5$-cross-polytopes. There exist lists of known finite polytopes for each kind (see \cite[Sections 12C,D,E]{arp}), and they are conjectured to be complete. 

For the first kind of the $6$-polytopes in \eqref{sixpol}, the underlying group $\Ga$ with generators $\rho_{0},\ldots,\rho_{5}$ has a presentation consisting of the standard Coxeter relations for the Schl\"afli type $\{3,3,3,4,3\}$ and the single extra relation
\begin{equation}  
\label{sixrel}
\left\{\begin{array}{rcccl}
(\rho_{1}\sigma\tau\sigma)^{s} & = & \epsilon, & \mbox{if} & k = 1,  \\
(\rho_{1}\sigma\tau)^{2s}               & = & \epsilon, & \mbox{if} & k = 2,
\end{array}\right.
\end{equation}
where 
\[ \sigma := \rho_{2}\rho_{3}\rho_{4}\rho_{3}\rho_{2}, \;\,
\tau := \rho_{5}\rho_{4}\rho_{3}\rho_{4}\rho_{5}. \]
This group $\Ga$ is known to be a finite C-group when ${\bf s}:=(2,0,0,0)$, $(2,2,0,0)$ or $(3,0,0,0)$, and the conjecture is that $\Ga$ is an infinite C-group for all other parameter vectors.

\begin{prob}
\label{unisix}
Complete the classification of the universal locally toroidal regular $6$-polytopes, that is, of the polytopes listed in \eqref{sixpol}.
\end{prob}

\medskip
\noindent
{\bf Locally toroidal chiral polytopes}
\smallskip

The facets and vertex-figures of a chiral polytope must be chiral or regular, and for a locally toroidal polytope the facets or vertex-figures must be toroidal as well. Thus locally toroidal chiral polytopes again can only exist for ranks $4$, $5$ or $6$, and must have the same Schl\"afli type as a locally toroidal regular polytope. In particular, up to duality, the rank $4$ types are those in \eqref{schla}, the only rank $5$ type is $\{3,4,3,4\}$, and the rank $6$ types possible are $\{3,3,3,4,3\}$, $\{3,3,4,3,3\}$ and $\{3,4,3,3,4\}$. Moreover, as there are no chiral toroids of rank $4$ or higher, a locally toroidal chiral polytope with chiral toroidal facets or vertex-figures must necessarily be of rank $4$. The corresponding universal polytope with these facets or vertex-figures then must be chiral as well. By contrast, there are examples of locally toroidal chiral polytopes whose facets and vertex-figures are regular (an example of rank $4$ has been described in \cite{cohupi}). In the case of regular facets and vertex-figures, the corresponding universal polytope with these facets and vertex-figures must necessarily be regular, and hence occurs among the polytopes discussed earlier. In particular, a locally toroidal universal polytope of rank $5$ or $6$ cannot be chiral; it must necessarily be regular.

This, then, limits the enumeration of the locally toroidal finite universal chiral polytopes with prescribed types of facets and vertex-figures to rank $4$ alone. Moreover, we may assume that at least one kind, facet or vertex-figure, is a chiral $3$-toroid of type $\{4,4\}_{(s,t)}$, 
$\{3,6\}_{(s,t)}$ or $\{6,3\}_{(s,t)}$ (with $t\neq 0$ or $s\neq t$). Bear in mind here that it does matter in which enantiomorphic form the facets and vertex-figures occur; thus there may be two (rather than one) essentially different amalgamation problems. The analysis again involves the seven possible (Schl\"afli) types $\{4,4,r\}$ with $r=3,4$, $\{6,3,r\}$ with $r=3,4,5,6$, and $\{3,6,3\}$.  

\begin{prob}
\label{chirtor}
Classify the universal locally toroidal chiral $4$-polytopes $\{\CP_1,\CP_2\}^{ch}$.  Here, one of $\CP_1$ and $\CP_2$ is a regular spherical or toroidal map and the other a chiral toroidal map, or both $\CP_1$ and $\CP_2$ are enantiomorphic forms of chiral toroidal maps.
\end{prob}

Problem~\ref{chirtor} translates into the enumeration of certain groups $\Ga$ defined in terms of generators and relations. These groups are quotients of the rotation (even) subgroup of the hyperbolic string Coxeter group on four nodes with branches labelled with the entries $p,q,r$ of the respective Schl\"afli symbol $\{p,q,r\}$. The three generators $\sigma_1,\sigma_2,\sigma_3$ for $\Ga$ satisfy the \emph{standard} relations
\[ {\sigma_1}^{p} = {\sigma_2}^{q} = {\sigma_3}^{r} =
{(\sigma_1\sigma_2)}^{2} = {(\sigma_2\sigma_3)}^{2} = {(\sigma_1\sigma_2\sigma_3)}^{2} = \epsilon ,\]
as well as one or two ``non-standard" relations determining the fine combinatorial and topological structure of the polytope.  For example, consider the universal chiral $4$-polytope $\{\{4,4\}_{(s,t)},\{4,3\}\}^{ch}$, whose facets are toroidal maps $\{4,4\}_{(s,t)}$ (the quotient of the square tessellation $\{4,4\}$ by the lattice spanned by the vectors $(s,t)$ and $(-t,s)$), and whose vertex-figures are $3$-cubes $\{4,3\}$.  Here there is just a single extra relation determined by $s$ and $t$, namely
\[ {(\sigma_1^{-1}\sigma_2)}^{s}  {(\sigma_1\sigma_2^{-1})}^{t} = \epsilon .\]
The main problem is to decide when $\Ga$ is finite. In this particular case it is conjectured that $\Ga$ is finite if and only $(s,t) = (1,2), (1,3), (1,4)$ or $(2,3)$ (up to interchanging $s$ and $t$); see \cite{swp} for supporting evidence. If the vertex-figure is also chiral, we obtain a similar relation from the vertex-figure, so then there are altogether four parameters $s,t,u,v$. Note that the corresponding universal regular polytope 
$\{\{4,4\}_{(s,t)},\{4,3\}\}$ (obtained when $t=0$ or $s=t$) is known to be finite if and only if $(s,t) = (2,0), (2,2)$ or $(3,0)$ (see \cite[Section~10B]{arp}). 

\medskip
\noindent
{\bf Other topological types}
\smallskip

Regular or chiral polytopes of other local topological types have not yet been systematically studied. On the other hand, many interesting examples are known, especially in rank $4$. Even if the topological type is preassigned only for one kind, facet or vertex-figure, very little is known about existing universal polytopes. The above Problems~\ref{resphvert} and \ref{chsphvert} address the interesting special case when the vertex-figure (say) is spherical.

In rank $4$, the facets and vertex-figures are maps on closed surfaces and the investigation can draw upon recent progress in the study of such maps. For orientable surfaces of genus at most $15$ and for non-orientable surfaces of genus at most $30$, the regular and chiral maps have recently been enumerated in \cite{conder}, complementing earlier work (for example, that in \cite{cm,ga,sh,wils78}) for maps of small genus; in particular, presentations of their groups have been found. (Recall that the genus $g$ of a closed surface $S$ is related to its Euler characteristic $\chi$ by $\chi =2-2g$ if $S$ is orientable, and $\chi = 2-g$ if $S$ is non-orientable.) 

The regular or chiral polytopes of rank $4$ naturally fall into four families specified by the following properties:\\[-.28in]
\begin{itemize}
\item facets and vertex-figures orientable,\\[-.28in]
\item facets and vertex-figures non-orientable,\\[-.28in]
\item facets orientable and vertex-figures non-orientable,\\[-.28in]
\item facets non-orientable and vertex-figures orientable.\\[-.25in]
\end{itemize}
Clearly, the last two families are equivalent under duality. The polytopes in each family can be further specified by the genera $g_1$ of their facets and $g_2$ of their vertex-figures. 

\begin{prob}
\label{othtop}
Let $g_{1},g_{2}$ be a pair of (small) non-negative integers. For each of the above four families, classify the finite universal regular or chiral $4$-polytopes $\{\CP_1,\CP_2\}$ and $\{\CP_1,\CP_2\}^{ch}$ with facets $\CP_1$ of genus $g_1$ and vertex-figures $\CP_2$ of genus $g_2$.
\end{prob}

The main interest is clearly in small genera $g_1$ and $g_2$. The locally toroidal $4$-polytopes belong to the first family and arise when $(g_{1},g_{2}) = (0,1)$, $(1,0)$ or $(1,1)$. 

\section{Realizations}
\label{realiz}

There are two main directions of research in the theory of regular and chiral polytopes:  the abstract, purely combinatorial aspect, and the geometric one of realizations. We now concentrate on the latter. In fact, much of the appeal of regular polytopes throughout their history has been their geometric symmetry. Here we restrict ourselves to realizations in Euclidean (or spherical) spaces, although some concepts generalize to other spaces (for example, hyperbolic spaces, or spaces over finite fields) as well. For general background see \cite[Chapter~5]{arp}.

Let $\CP$ be an abstract $n$-polytope, and let $\CP_j$ denote its set of $j$-faces.  Following \cite[Section~5A]{arp}, a {\em realization\/} of $\CP$ is a mapping $\beta\colon\,\CP_{0}\rightarrow \BE^{d}$ of the vertex-set $\CP_{0}$ into some Euclidean space $\BE^d$.  Define $\beta_0 := \beta$ and $V_{0} := V := \CP_{0}\beta$, and write $2^X$ for the family of subsets of the set $X$. The realization $\beta$ recursively induces surjections $\beta_{j}\colon\, \CP_{j} \rightarrow V_{j}$, for $j=1,2,3$, with $V_{j}\subset 2^{V_{j-1}}$ consisting of the elements
\[ F\beta_{j} := \{G\beta_{j-1} \mid G\in \CP_{j-1} \mbox{ and } G\leq F\} \]
for $F\in\CP_j$; further, $\beta_{-1}$ is given by $F_{-1}\beta_{-1} := \emptyset$. Even
though each $\beta_j$ is determined by $\beta$, it is helpful to think of the realization as
given by all the $\beta_j$. The realization is $d$-\emph{dimensional} if $\BE^d$ is the affine hull of $V$. We call the realization $\beta$ {\em faithful\/} if each $\beta_j$ is a
bijection. For a faithful realization, each $j$-face of $\CP$ with $j\geq 1$ must be uniquely determined by the $(j-1)$-faces which belong to it; if this purely combinatorial condition fails, $\CP$ does not admit a faithful realization. Note that a realization of $\CP$ determines a realization of each of its faces or co-faces. 

Our main interest is in discrete and faithful realizations. In this case, the vertices, edges, $2$-faces, etc., of $\CP$ are in one-to-one correspondence with certain points, line segments, simple (finite or infinite) polygons, etc., in $\BE^d$, and it is safe to identify a face of $\CP$ and its image in $\BE^d$. The resulting family of points, line segments, polygons, etc., is denoted by $P$, and it is understood that $P$ inherits the partial ordering of $\CP$; when convenient $P$ will be identified with $\CP$. The \emph{symmetry group} $\CG(P)$ of a $d$-dimensional realization $P$ of $\CP$ is the group of all isometries of $\BE^d$ that maps $P$ to itself. 

For a faithful realization of an abstract regular $n$-polytope $\CP$ (as a geometrically regular polytope) we have two ingredients. First, we need a faithful representation of its automorphism group $\Ga(\CP)$ as a group $\CG$ of isometries on some Euclidean space $\BE^d$; this group $\CG$ is the symmetry group of the realization of $\CP$ and is generated by the images $R_{0},\ldots,R_{n-1}$ in $\CG$ of the distinguished generators $\rho_{0},\ldots,\rho_{n-1}$ of $\Ga(\CP)$.  The generators $R_j$ of $\CG$ are reflections in subspaces, their \emph{mirrors}, of $\BE^d$. The subspace $W$ fixed by $R_{1},\ldots,R_{n-1}$ is called the \emph{Wythoff space} of the realization. The realization of $\CP$ associated with $\CG$ and its generators $R_{j}$ then arises from some choice of \emph{initial vertex} $v$ in $W$; its vertex-set is $V := v\CG$, the orbit of $v$ under $\CG$. The actual realization $P$ of $\CP$ is obtained by \emph{Wythoff's construction}:  the faces in its base flag are $F_0 := v$, and, for $j\geq 1$, 
\[  F_j := F_{j-1} \langle R_0,\ldots,R_{j-1} \rangle,  \]
and its $j$-faces are the $F_jR$ with $R\in\CG$, with the order relation given by iterated membership.  The realization $P$ then has a flag-transitive symmetry group $\CG(P):=\CG$ and is often called a \emph{regular geometric polytope\/}.

There is a similar realization theory for chiral polytopes. A realization $P$ of an abstract polytope $\CP$ is called \emph{chiral} (or a \emph{chiral geometric polytope\/}) if $P$ has two orbits of flags under its symmetry group $\CG(P)$, with adjacent flags lying in different orbits (see \cite{rcld,esch2}). Here the original polytope $\CP$ must be (combinatorially) regular or chiral. Chiral realizations are derived by a variant of Wythoff's construction, applied to a suitable representation $\CG = \langle S_1,\ldots,S_{n-1} \rangle$ of the underlying combinatorial group $\Ga := \langle \sigma_{1},\ldots,\sigma_{n-1} \rangle$ (with $S_{i}$ the image of $\sigma_i$); the latter is $\Ga(\CP)$ or $\Ga^{+}(\CP)$ according as the abstract polytope $\CP$ is chiral or regular. The \emph{Wythoff space} now is the fixed set of the subgroup $\CG_0 := \langle S_2,\ldots,S_{n-1} \rangle$. An abstract regular polytope may have chiral realizations, though not necessarily faithful ones; it is an interesting open question whether it could actually have faithful chiral realizations.  

There are two main directions of research in the realization theory of abstract polytopes, and most enumeration projects of symmetric polytopes that have been undertaken in the past follow one of these approaches. 

The first, for which a fairly complete general theory exists for regular polytopes (at least in the finite case), studies the space of all (regular) realizations of a given abstract regular polytope $\CP$; for finite polytopes, this realization space has the structure of a closed convex cone, the \emph{realization cone} of $\CP$ (see \cite[Chapter 5]{arp} and \cite{hemo,mcreal,mcmon}), and its fine structure is determined by the family of irreducible orthogonal representations of the automorphism group $\Ga(\CP)$. Here, much less is known about the space of chiral realizations of chiral or regular abstract polytopes. 

\begin{prob}
\label{devchi}
Develop the details of the realization theory for chiral polytopes.
\end{prob}

The second, traditional approach asks for the classification of the realizations of all polytopes in a Euclidean space of given dimension, and it is usual to impose the condition that the realization be discrete and faithful.  

\medskip
\noindent
{\bf Realization cones and real representations\/}
\vskip.02in

Each realization of a finite regular polytope $\CP$ in a Euclidean space is uniquely determined by its \emph{diagonal vector\/}, whose components are the squared lengths of the diagonals (pairs of vertices) in the diagonal classes of $\CP$ modulo $\Ga(\CP)$.  A realization can be identified with its diagonal vector, and then the realization cone of $\CP$ simply consists of all possible diagonal vectors of realizations of $\CP$. Each orthogonal representation $\CG$ of $\Ga(\CP)$ yields a (possibly degenerate) realization of $\CP$. The degree of freedom for deriving a realization from $\CG$ is measured by the dimension of the \emph{Wythoff space\/}.  The sum of two diagonal vectors in the cone corresponds to the {\it blend\/} of the two realizations, and this, in turn, to the sum of the two representations of $\Ga(\CP)$. The \emph{pure\/} (non-blended) realizations of $\CP$ are determined by the irreducible representations of $\Ga(\CP)$; they correspond to the extreme rays of the cone.  There are numerical relationships that involve the dimensions of the Wythoff spaces and the degrees of the irreducible representations of $\Ga(\CP)$ (see \cite[Section 5B]{arp}). Some relationships have recently been explained using character theory (see \cite{hemo}). 

The fine structure of the realization cone is only known for a small number of polytopes, including the regular convex polytopes, with the exception of the $120$-cell $\{5,3,3\}$ and $600$-cell $\{3,3,5\}$, and the regular toroids of rank $3$ (see \cite[Section 5B]{arp} and \cite{bust,mowe,monw5}).  

\begin{prob}
\label{realtor}
Describe the realization cone for the regular toroids of any rank.
\end{prob}

The groups of the regular toroids have large abelian subgroups of relatively small index, so the representation theory may still be manageable in this case. 

\begin{prob}
\label{realcharac}
Can a finite string C-group $\Gamma$ have an irreducible representation $T:\Gamma \rightarrow \rm{GL}_{n}(\BC)$ with real character, but which is not similar to a real representation? (This means that the Frobenius-Schur indicator of $T$ is $-1$.)
\end{prob}

In case of a positive answer to Problem~\ref{realcharac}, the representation $T\oplus\overline{T}$ associated with the pair $T,\overline{T}$ of complex conjugate representations, may provide an interesting realization of the regular polytope $\CP$ corresponding to $\Gamma$ (see \cite{hemo}). 

Our next problem is of general interest, independent of realization theory.

\begin{prob}
\label{intcharac}
Suppose $\Gamma = \langle \rho_{0},\ldots,\rho_{n-1} \rangle$ is a finite group generated by involutions with string diagram (that is, the $\rho_{i}$ satisfy \eqref{standrel}). Can one determine the validity of the intersection condition \eqref{intprop} for $\Gamma$ from its character table?  Or, can the validity be determined from some other readily computed invariants of the group algebra?
\end{prob}

\medskip
\noindent
{\bf Realizations in a space of given dimension\/}
\vskip.02in

There has been considerable progress on realizations of regular or chiral polytopes in small dimensions since the publication of \cite{arp}; for a brief survey see \cite{rcld}. The regular polyhedra in $\BE^3$, also known as \emph{Gr\"unbaum-Dress polyhedra}, were already enumerated well over twenty-five years ago in \cite{dress, gruen} (see \cite[Section 7E]{arp} for an alternative approach), but the chiral polyhedra in $\BE^3$ were only described recently in \cite{esch2}. The regular polytopes of \emph{full rank}, which comprise the finite regular polytopes of rank $n$ in $\BE^n$ as well as the infinite discrete regular polytopes of rank 
$n+1$ in $\BE^n$, were classified by McMullen~\cite{full}. Moreover, the finite regular polyhedra in $\BE^4$ were completely enumerated in \cite{fourdim}.

There are many other classes of more or less symmetrical polyhedra or polytopes in Euclidean spaces, whose symmetry groups have transitivity properties which are weaker than flag-transitivity (see \cite{mar} for an overview of classical efforts in this direction). Most enumeration problems are open even for ordinary $3$-space.

\begin{prob}
\label{fulltreu}
Enumerate the finite or infinite discrete polyhedra in $\BE^3$ which are geometrically fully transitive (that is, their symmetry group acts transitively on the vertices, edges, and the $2$-faces).
\end{prob}

Here even the case of finite polyhedra is not settled; see \cite{far} for partial results on the enumeration of finite fully transitive polyhedra in $\BE^3$.  Note that the regular and chiral polyhedra in $\BE^3$ are fully transitive, so any classification attempt, which necessarily must include these polyhedra, is likely to be an ambitious project.

Another relaxation of symmetry leads to the uniform polyhedra and polytopes. They have attracted a lot of attention and are the subject of the forthcoming monograph~\cite{jo}.
Recall that a polyhedron in $\BE^3$ is \emph{uniform} if its $2$-faces are regular polygons and its symmetry group is transitive on the vertices. The uniform convex polyhedra in $\BE^3$ consist of the five Platonic solids, the thirteen polyhedra known as Archimedean polyhedra, and two infinite classes of prisms and antiprisms. The finite uniform polyhedra with planar faces were first described in \cite{clm}, but the completeness of the list was only established (independently) in \cite{ski} and \cite{sop}. 

\begin{prob}
\label{fulltreu}
Enumerate the finite or infinite discrete, uniform polyhedra in $\BE^3$.
\end{prob}

In the finite case it remains to describe the uniform polyhedra with skew faces. For infinite discrete polyhedra very little is known. For non-faithful uniform realizations of abstract polyhedra in $\BE^3$ the reader is also referred to \cite{grnup}.  

Another interesting project, as yet unexplored, is the study of \emph{abstract uniform polytopes\/}.  As a first step this would require to phrase a sensible definition of the term ``uniform", and then develop the abstract theory and the geometric realization theory of such polytopes.

We conclude this section with a general remark about geometric polyhedra. A key element in 
Gr\"unbaum's~\cite{gruen} skeletal approach to regular polyhedra in $\BE^3$ was the idea to restore the symmetry in the definition of regularity by allowing the $2$-faces, not only the vertex-figures (as in \cite{skew}), to be skew polygons. Another digression from the traditional approach was to admit infinite, zigzag or helical $2$-faces. These radically new ideas provided us with new classes of polyhedra that have been studied in detail in the context of symmety but have otherwise not attracted the attention they deserve. It seems to us that a fascinating world of skeletal polyhedra awaits to be discovered. Suffice it here to end with a question that points to some of the possibilities:\ what can be said about the general properties of helix-faced polyhedra? (Most regular, and all chiral, infinite-faced polyhedra in $\BE^3$ have helical $2$-faces.)
 
\section{Other problems}
\label{othpro}

In this section we collect the remaining problems which do not fit under the headings of the previous sections.

\medskip
\noindent
{\bf Medial layer graphs\/}
\vskip.02in

\begin{figure}[hbt]
\begin{center}
\begin{picture}(150,220)
\put(-70,10)
{\includegraphics[width=4in]{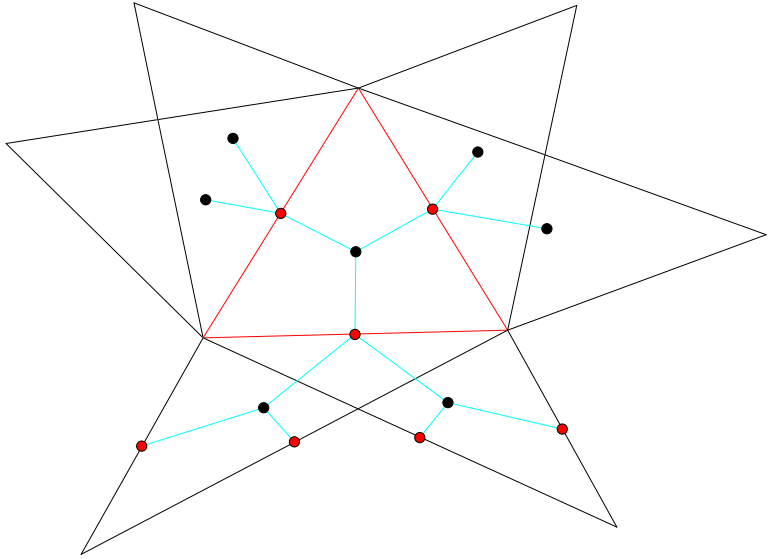}}
\end{picture}
\caption{\em The medial layer graph of a polytope of type $\{3,q,3\}$.}
\label{figtwo}
\end{center}
\end{figure}

When $\CP$ is a finite self-dual regular (or chiral) $4$-polytope of type $\{3,q,3\}$, the faces of ranks $1$ and $2$ can be thought of as the vertices of a finite bipartite, trivalent, $3$-transitive (or $2$-transitive, respectively) graph $\CG$ called the \emph{medial layer graph\/} of $\CP$; in $\CG$, two vertices are joined by an edge if and only if the corresponding faces in $\CP$ are incident. Recall that a trivalent graph is $t$-\emph{transitive} for some $t$ if its automorphism group is transitive on $t$-arcs, but not on $(t+1)$-arcs. The construction of the graph $\CG$ from the polytope $\CP$ is illustrated in Figure~\ref{figtwo}, where the vertices of $\CG$ in the two classes are represented by the midpoints of edges or centers of $2$-faces of $\CP$, respectively. 

Given such a graph $\CG$, a reverse construction yields a poset $\CP(\CG)$ of rank $4$ with many properties characteristic for polytopes (see \cite{mow}).

\begin{prob}
\label{pg}
Given a finite, bipartite, trivalent, $3$-transitive or $2$-transitive graph $\CG$, when is the poset $\CP(\CG)$ a polytope?
\end{prob}

The solution would amount to giving properties intrinsic to the graph $\CG$, which would guarantee the reconstruction of a polytope.

When the polytope $\CP$ is not self-dual, its medial layer graph $\CG$ is not symmetric (arc-transitive) but only semisymmetric (edge- but not vertex-transitive); see \cite{mpsw} for interesting semisymmetric graphs obtained in this way.

\medskip
\noindent
{\bf Petrie schemes\/}
\vskip.02in

Petrie schemes are combinatorial substructures in polytopes that generalize Petrie polygons (see \cite{gow1,gow2}). Let $\CP$ be an abstract $n$-polytope, let $\CF$ be its set of flags, and let $\gamma_{i}$ ($i=0,\ldots, n-1$) denote the involutory bijection of $\CF$ which maps \emph{each} flag onto its $i$-adjacent flag. (The $\gamma_i$'s are not automorphisms of $\CP$.) Consider bijections of $\CF$ of the form $\gamma := \gamma_{i_{1}}\gamma_{i_{2}}\ldots\gamma_{i_{n}}$, where $i_{1},i_{2},\ldots,i_{n}$ is a permutation of $0,1,\ldots,n-1$. From any such product $\gamma$ and any flag $\Psi$ of $\CP$ we obtain a ($2$-sided) infinite sequence of flags
\[(\ldots,\Psi\gamma^{-2},\Psi\gamma^{-1},\Psi,\Psi\gamma,\Psi\gamma^2,\ldots). \] 
A \emph{Petrie scheme\/} of $\CP$ is the shortest possible representation of any such sequence. More precisely, if the sequence contains repeating cycles of flags, then the Petrie scheme is the shortest possible cycle presentation of that sequence; otherwise, the Petrie scheme is the sequence itself. A Petrie scheme is \emph{acoptic} if each proper face of $\CP$ appears at most once in the flags of the Petrie scheme (that is, the Petrie scheme ``has no self-intersections").  A polytope can have both acoptic Petrie schemes and non-acoptic Petrie schemes. 

It is an interesting problem to determine which polytopes have only acoptic Petrie schemes. This is true for the regular convex polytopes, the regular Euclidean tessellations, and the regular polyhedra in $\BE^3$ (see \cite{gow1,gow2}). Here we collect two problems addressing this property for other classes of polytopes, including polytopes with less symmetry. 

\begin{prob}
\label{acopcl}
Which of the following types of convex polytopes or star-polytopes have the property that all their Petrie schemes are acoptic?\\[-.31in]
\begin{itemize}
\item Cartesian products of cubes and crosspolytopes. \\[-.31in]
\item Cartesian products of cubes or crosspolytopes and $4k$-gons, $k\geq 1$.\\[-.31in]
\item Uniform convex polytopes.\\[-.31in]
\item Regular star-polytopes.\\[-.31in]
\end{itemize}
\end{prob}

The answer is conjectured to be affirmative at least for the first two and the last types (see 
\cite{gow2}).

\begin{prob}
\label{acopuni}
Which universal regular polytopes $\{\CP_1,\CP_2\}$ have the property that all their Petrie schemes are acoptic?
\end{prob}

Not all regular polytopes have only acoptic Petrie schemes, so at the minimum we must require here that the facet $\CP_1$ and vertex-figure $\CP_2$ (of rank $n$) possess no Petrie schemes that self-intersect within $n+1$ steps.

\medskip
\noindent
{\bf Specific groups}
\smallskip

Of great interest, particularly to researchers new to the field, are the following two web-based atlases:\  ``The Atlas of Small Regular Polytopes" in~\cite{hatlas}, and ``An atlas of abstract regular polytopes for small groups" in~\cite{leem}. 

Next we mention a problem that arose in the context of creating the atlas in \cite{hatlas}. When trying to determine which groups of a given (small) order are automorphism groups of polytopes, the orders $2^k$ or $2^{k}p$ proved to be more difficult than others.

\begin{prob}
\label{grorder}
Characterize the groups of orders $2^k$ or $2^{k}p$, with $k$ a positive integer and $p$ an odd prime, which are automorphism groups of regular or chiral polytopes.
\end{prob}

Abstract polytope theory has interesting connections with incidence geometries, diagram geometries, and buildings (see \cite{bup}). It is not the appropriate place here to elaborate on this interplay in detail. Suffice it to mention the following problem, which is a polytope version of a more general theme that has driven the development of the theory of diagram geometries.

\begin{prob}
\label{spor}
Find regular, chiral, or other polytopes whose automorphism groups are sporadic simple groups. For example, is there a polytope whose automorphism group is the Monster group?
\end{prob}

A similar problem can be posed for arbitrary finite simple groups. The question which finite simple groups occur as automorphism groups of regular maps, is already contained in the Kourovka Notebook~\cite{kou} as Problem 7.30 (see also \cite{ned1}). 

The following problem (asked by Michael Hartley) deals with alternating groups; for solutions in rank $3$ see the recent papers \cite{ned2,pel}.

\begin{prob}
\label{altgr}
Find regular, chiral, or other polytopes whose automorphism groups are alternating groups $A_n$. In particular, given a rank $r$, for which $n$ does $A_n$ occur as the automorphism group of a regular or chiral polytope of rank $r$?
\end{prob}

Call an $n$-polytope $\CP$ \emph{neighborly\/} if any two of its vertices are joined by an edge. There are generalizations of this property to $k$-neighborliness, requiring that any $k$ vertices of $\CP$ are the vertices of a $(k-1)$-face. Here we restrict ourselves to $k=2$, as the following problem is already open in this case. 

\begin{prob}
\label{spor}
Characterize the finite neighborly, regular or chiral polytopes.
\end{prob}

The automorphism group of a finite neighborly, regular or chiral polytope $\CP$ must necessarily act $2$-transitively on the vertices of $\CP$. Thus the characterization problem is closely related to the enumeration of $2$-transitive permutation groups.

\medskip
\noindent
{\bf Quotient polytopes}
\smallskip

For a string C-group $\Ga = \langle\rho_{0},\ldots,\rho_{n-1}\rangle$ we use the following notation. For $i=0,1,\ldots,n-1$, define $\Ga_{<i}:=\langle\rho_{j}\mid j<i\rangle$ and $\Ga_{>i}:=\langle\rho_{j}\mid j>i\rangle$, and recall that $\Ga_{i}:=\langle\rho_{j}\mid j\neq i\rangle$. Moreover, for a subgroup $\Sigma$ of $\Ga$ and $\varphi\in\Ga$ we write 
$\Sigma^{\varphi}:=\varphi^{-1}\Sigma\varphi$.

The \emph{quotient} of an $n$-polytope $\CP$ by a subgroup $\Sigma$ of its automorphism group $\Ga(\CP)$ is the set of orbits $\CP/\Sigma$ of faces of $\CP$ under the action of $\Sigma$, with two orbits incident if they contain incident faces of $\CP$ (see \cite[Section 2D]{arp} and \cite{ha97}). A \emph{quotient polytope} is a quotient which is again a polytope.

Any abstract $n$-polytope $\CQ$ may be constructed as a quotient of a regular $n$-polytope $\CP$ by a semisparse subgroup of its automorphism group $\Ga(\CP)$ (see \cite{ha97}). Recall that a subgroup $\Sigma$ of a string C-group $\Ga (=\Ga(\CP))$ is \emph{semisparse} if and only if the following three properties are satisfied:\\[-.26in]
\begin{itemize}
\item $\Sigma^{\varphi} \,\cap\, \Ga_{i}\rho_{i} = \emptyset$\  for all $\varphi\in\Ga$ and 
$0 \leq i \leq n-1$,\\[-.26in]
\item $\Ga_{>i} \,\cap\, \Ga_{<j}\Sigma^{\varphi} \,\subseteq\,
\Gamma_{k} (\Ga_{>i} \cap \Ga_{<j}) \Sigma^{\varphi}$\  for all $\varphi\in\Ga$ and 
$0 \leq i < k < j \leq n-1$,\\[-.26in]
\item  $\Ga_{>i} \,\cap\, \Ga_{<i+1}\Sigma \,\subseteq\,  \Sigma$\  for all $0 \leq i < n-1$.
\\[-.26in]
\end{itemize} 
It is generally difficult to tell if a subgroup $\Sigma$ of a string C-group is semisparse. By \cite[Thm. 2.4]{hasemi}, if $\Sigma$ is a semisparse subgroup of $\Ga$, then, for all $\varphi\in\Ga$, the subgroups $\Sigma^{\varphi} \cap \Ga_{0}$ and $\Sigma^{\varphi} \cap \Ga_{n-1}$ are semisparse in $\Ga_0$ and $\Ga_{n-1}$, respectively, and
\begin{equation}
\label{semisparse}
\Sigma^{\varphi} \,\cap\, \Ga_{0}\Ga_{n-1} 
\,=\, (\Sigma^{\varphi} \cap \Ga_{0})(\Sigma^{\varphi} \cap \Ga_{n-1})
\quad (\varphi\in\Ga).
\end{equation}
The vertex-figures and facets of the corresponding quotient polytope $\CP/\Sigma$ then are given by the quotients $\CP^{0}/(\Sigma^{\varphi} \cap \Ga_{0})$ and 
$\CP^{n-1}/(\Sigma^{\varphi} \cap \Ga_{n-1})$ for the various elements $\varphi\in\Ga$, where here $\CP^0$ and $\CP^{n-1}$ denote the vertex-figures and facets of $\CP$ (associated with $\Ga_{0}$ and $\Ga_{n-1}$), respectively.

Our last problem calls for the proof of sufficiency of the conditions in \eqref{semisparse}; see also \cite{hasemi} for partial results in this direction.

\begin{prob}
\label{convsemi}
Let $\Sigma$ be a subgroup of $\Ga$. Is it true that whenever $\Sigma^{\varphi} \cap \Ga_{0}$ and $\Sigma^{\varphi} \cap \Ga_{n-1}$, respectively, are semisparse subgroups of $\Ga_0$ and $\Ga_{n-1}$ for all $\varphi\in\Ga$, and condition \eqref{semisparse} holds, then $\Sigma$ is a semisparse subgroup of $\Ga$?
\end{prob}

\noindent
{\bf Acknowledgment:}  We would like to thank all participants of the two meetings at BIRS and Calgary for their contributions to the program. We are especially grateful to those who shared problems with us at the problem session. Last but not least we would like to thank Ted Bisztriczky for generously providing us with the opportunity to follow the BIRS Workshop by the Polytopes Day at the University of Calgary.

\vskip.2in
\noindent
{\em Addresses\/}:\\[.1in]
{\em Egon Schulte, Northeastern University, Boston, MA 02115, USA. }\\[.03in] 
{\em Asia Ivi\'{c} Weiss, York University, Toronto, ON, Canada M3J 1P3.}

\vskip.2in
\noindent
{\em Email addresses\/}:\\[.1in]  
schulte@neu.edu,\  weiss@mathstat.yorku.ca\\[.1in]
\end{document}